\documentclass[12 pt]{article}
\fontsize{13.5}{42}
\usepackage{amsmath,amssymb,braket}
\usepackage[english]{babel}
\usepackage{lipsum}
\usepackage{hyperref}
\date{}
\begin{document}
\title{\textbf{SOME OBSERVATIONS ON A RESULT BY BIA\L{}INICKI-BIRULA AND \.{Z}ELAZKO}}
\author{\sc Gabriele Gull\`a}
\maketitle
\begin{abstract} 
We are interested in investigating some definitions and assumptions stated in \cite{BBZ}, in particular the notions of \textit{measurability} and \textit{atomicity} that the two authors used in order to give a representation for multiplicative linear functionals (m.l.f. so on) over the cartesian product of algebras.\\
We want to briefly analyze the result under other strong proper axioms with respect to $ZF(C)$.\\
Some little historical and notational remarks will be given.
\end{abstract}
\textbf{Key Words}: Measurable cardinals, Real-valued measurable cardinals, 2-measurable cardinals, multiplicative linear functionals.
{\let\thefootnote\relax}
\section{Introduction}
As in \cite{BBZ}, we fix for the rest of the paper the following: let $\mathbb{K}$ be a field, $T$ an arbitrary abstract set, $\mathcal{R}_t$, $t\in T$, a family of unitary algebras over $\mathbb{K}$; let us define
$$\mathcal{R}=\prod_{t\in T}\mathcal{R}_t$$ 
The elements of $\mathcal{R}$ can be seen as functions defined over $T$ which maps every $t\in T$ to an element $x(t)\in \mathcal{R}_t$.\\
If $f_{t_0}$ is a m.l.f. on $\mathcal{R}_{t_0}$, then\\\\
(1)$\qquad\qquad\qquad\qquad f(x)=f_{t_0}(x(t_0))$\\\\
is a m.l.f. defined on $\mathcal{R}$.\\
\\Bia\l{}inicki-Birula and \.{Z}elazko showed in \cite{BBZ} the following result, that we present in the same way of the paper:\\
\\\textit{Theorem 1.1}\footnote{\textit{Fundamental theorem}, page 2 of \cite{BBZ}}: \textit{If $|\mathbb{K}|<\aleph_I$ (\textit{the first aleph of measure $\neq 0$}), then every m.l.f. defined over $\mathcal{R}$ is given by
$(1)$ iff every $(0 - 1)$ $\sigma$-measure defined over the field of all subset of $T$ is the \textit{atomic one}, i.e. $|T|<\aleph_I$}.\\
\\\textbf{Remark 1.1}\\
The two authors proved that, if $|\mathbb{K}|\ge|T|$, every m.l.f. can be expressed by (1)\footnote{Theorem 3 in \cite{BBZ}} so, in what follows, we assume that $|\mathbb{K}|<|T|$.\\
\\I want to start by investigating both historically and mathematically the notion of cardinal measurability used as hypothesis in \cite{BBZ}.
\section{Measurability properties}
The authors gave the following definition (again, stated here almost in the same way of the paper):\\
\\\textit{Definition 2.1}\footnote{First footnote of page 2 of \cite{BBZ}}: \textit{A measure is called  $(0 - 1)$ if it assume only two values: 0 and 1.
A cardinal $\kappa$ is said to be of measure $\neq 0$ if there exists a $(0 - 1)$ $\sigma$-measure $\mu$ defined on the field of all subset of set of power $\kappa$ which vanishes on the points.}\\
\\Here ``vanishing on the points'' means that the measure is nontrivial, i.e. $\mu(\left\{x\right\})=0$ for every singleton $\left\{x\right\}$. Being non-trivial, for a measure, is a property almost always implicitly required.\\Following, for example, Scott\cite{S}, Choquet\cite{C} or Fremlin\cite{F}, we call the above mentioned kind of cardinals \textit{2-(valued)measurable cardinals} and a $(0 - 1)$ measure is a \textit{two-valued measure}.\\
By the definition it follows that every cardinal greater than the first\\ 2-measurable is still 2-measurable, and so this (anyway large) cardinal notion, from a set-theoretic point of view, is not - let say - \textit{directly} interesting as others since the first 2-measurable it just divides the class of cardinals in two segments.\\
\\The first notion of measurability for cardinals was introduced by Banach-Kuratowski\cite{BK} and Banach\cite{B} between 1929 and 1930 in order to generalize and solve the ``\textit{probl\`eme de la mesure}'' introduced by Lebesgue's works.\\
The modern definition of real-valued measurable cardinal is the following:\\
\\\textbf{Definition 2.2}\\An uncountable cardinal $\kappa$ is \textit{real-valued measurable} if there exists on it a nontrivial probabilistic $\kappa$-additive measure.\\
\\In 1930, Ulam gave in \cite{U} the definition of measurable cardinal:\\
\\\textbf{Definition 2.3}\\An uncountable cardinal $\kappa$ is \textit{measurable} if there exists on it a $\kappa$-complete nonprincipal ultrafilter.\\
It turns out that every measurable cardinal is \textit{(strongly) inaccessible}\footnote{see \cite{J} or \cite{K}.}.\\
\\The relations amongst these notions of measurable-like cardinals are the following:\\
(a) \centerline{$\kappa$ measurable $\Rightarrow$ $\kappa$ 2-measurable}\\
(b) \centerline{$\kappa$ measurable $\Rightarrow$ $\kappa$ real-valued measurable}\\
(c)\footnote{The converse holds under $CH$.} \centerline{$\kappa$ is the first 2-measurable $\Rightarrow$ $\kappa$ is the first measurable}\\
(d) \centerline{$\kappa$ is the first 2-measurable $\Rightarrow$ $\kappa$ real-valued measurable}\\
(e)\footnote{The first real-valued measurable can be smaller than the first 2-measurable because, under some assumptions, there are real valued measurable less or equal than $\mathfrak{c}$ and the latter, under $CH$, is not measurable.} \centerline{$\mathfrak{c}<\kappa$ real-valued measurable $\Rightarrow$ $\kappa$ measurable}\\
\\
\\For (a): if $\kappa$ carries a non principal $\kappa$-complete ultrafilter $\mathfrak{U}$ (which is of course $\sigma$-complete) then the measure
$$\mu_{\mathfrak{U}}(X)=\begin{cases}1\qquad X\in \mathfrak{U}\\0\qquad X\notin \mathfrak{U}\end{cases}$$
is a two-valued $\sigma$-measure.\\
\\For (b): trivial, since every two-valued measure is a probabilistic one. \\
\\For (c) and (d): if $\kappa$ carries a two-valued non trivial $\sigma$-measure $\mu$, the set
$$\mathfrak{U}_{\mu}=\left\{X\subset \kappa: \mu(X)=1\right\}$$
is a nonprincipal $\sigma$-complete ultrafilter on $\kappa$: $\sigma$-completeness comes from the $\sigma$-additivity of the measure and non principality from the nontriviality.\\
Let $\kappa$ be the first cardinal with such a property; it follows that the ultrafilter $\mathfrak{U}_{\mu}$ is $\kappa$-complete too\footnote{see Ulam\cite{U}.}, and so $\kappa$ is (the first, from (a)) measurable and (from (b)) real-valued measurable.\\
\\We will see (e) in the next section.\\
\\\textbf{Remark 2.1}\\
We want to stress the fact that here we are working with respect to the theory $ZFC$, adding the proper large cardinal axioms $RVM$ (``There exists a real-valued measurable cardinal'') and $M$ (``There exists a measurable cardinal'')\footnote{From the consistency point of view Solovay\cite{So} has proven that $Con(ZFC+M)\Leftrightarrow Con(ZFC+RVM)\Leftrightarrow Con(ZFM)$.}.\\
\\It follows that the first 2-measurable $\kappa$ is both measurable and real-valued measurable, and so, in Theorem 1.1, the field $\mathbb{K}$ can be actually the reals (a relevant case) and $T$ and $\mathbb{K}$ can have real-valued measurable cardinalities.
\section{Atomicity conditions}
If $\mu$ is a measure on a set $X$ such that there exists an $A\subset X$ with:\\1) $\mu(A)>0$ and\\
2) $\forall B\subset A\enspace (\mu(B)=0 \vee \mu(B)=\mu(A))$\\
then $A$ is called \textit{atom} for the measure, which then is called \textit{atomic}. A measure without atoms is called \textit{non atomic} or \textit{atomless}.\\
\\About (e) of the previous section: if $\kappa$ carries an atomless measure then $\kappa$ is $\le \mathfrak{c}$ (see Lemma 27.5 of \cite{J} for details); if, \textit{au contraire}, the measure $\mu$ has an atom $A$, the set
$$\mathfrak{U}_{\mu}=\left\{X\subset \kappa: \mu(X\cap A)=\mu(A)\right\}$$
is a complete non principal ultrafilter, then we can define a measure $\nu$ by $\nu(X) = \mu(X\cap A)/\mu(A)$. Then $\nu$ is a two-valued, $\kappa$-additive measure and $\kappa$ is measurable.\\
\\Now let us concentrate on the second part of Theorem 1.1: ``\textit{...every $(0 - 1)$ $\sigma$-measure defined over the field of all subset of $T$ is the atomic one, i.e. $|T|<\aleph_I$}''.\\\\
Let $\mu$ be a two-valued measure on $X$ and let $A$ be such that $\mu(A)>0$. Then, $\mu(A)=1$. For any subset $B$ of $A$, either $\mu(B)=0$, or $\mu(B)=1$. So, $A$ is an atom and $\mu$ is atomic\footnote{In particular in 1922 Sierpinski\cite{S} proved that non-atomic measures have continuum many values.}.\\
It follows that the  the second member of the equivalence stated in Theorem 1.1 is always vacuously true for every set carrying a two-valued measure.\\
Moreover the authors, proving the theorem, wrote: ``\textit{We now suppose that on the set $T$ there exists a non-atomic $(0 - 1)$ $\sigma$-measure $\mu$...}'' which then turns to be inconsistent.\\
\\So it arises the question about the notion of atomicity used by the authors.\\By definition of 2-measurability it follows that, if a set $X$ has cardinality less that $\aleph_I$ (or $\aleph_{\mu}$, following the conventional modern notation for the first measurable), then every two-valued $\sigma$-measure does not vanish on the points, meaning it is trivial. So there exists a point $x$ such that $\mu(\left\{x\right\})=1$.\\
And that is in fact what is meant when they write: ``\textit{...$\tilde{\mu}$ must be atomic (...) and there exists a point $p$ (...) such that $\tilde{\mu}(\left\{p\right\})=1$}''\footnote{Page 4 of \cite{BBZ}, proof of Theorem 2.}.\\
\\It turns out that Theorem 1.1 can be simply stated as follows\\
\textit{Theorem 1.1*}: If a field $\mathbb{K}$ has cardinality smaller than the first measurable cardinal, then every m.l.f. defined over $\mathcal{R}$ is given by $(1)$ iff the set $T$ has cardinality smaller than the first measurable cardinal.
\section{Constructibility and measurability}
In 1961 Dana Scott\cite{S} proved the following result:\\\\
\textbf{Theorem 4.1} (Scott, 1961)\\
The theory $ZFC$+V=L+``There is a 2-measurable cardinal'' is inconsistent.\\
\\It follows that, under V=L, there are not 2-measurable cardinals and so:\\
\\\textbf{Corollary 4.2}\\
$ZFC$+V=L $\vdash$ ``Every m.l.f. defined over $\mathcal{R}$ is given by $(1)$ for every field $\mathbb{K}$ and every set $T$.''\\
\\\textbf{Remark 4.1}\\
If we replace the Axiom of Choice with the \textit{Axiom of Determinacy} (they exclude each other) we obtain the theory $ZF+AD$ which proves that $\aleph_1$ is measurable\footnote{see \cite{J2}.}.\\
Under $AD$ we have that $\aleph_1<\mathfrak{c}$, and so, in this context, $\mathbb{K}$ must be finite and $T$ can be at most countable.
\subsection{Appendix: historical and notational remarks}
We want to conclude with a folklore note; in fact it is interesting the notation used in \cite{BBZ} for the sets: in order to represent a collection
$$\left\{x\in X: P(x)\right\}$$
it has been used the following construct
$$\underset {x\in X}E\enspace P(x)$$
We are not aware of other contexts which used this symbology, but it seems to be a residual heritage of french culture influence which was - with the german one - quite relevant in eastern Europe before the spreading of english language after the II WW (which in fact was used as the language for \cite{BBZ}).\\
The symbol $E$, in fact, stands for the french \textit{ensemble}, underlined by the kind of elements and followed by the predicate which describes them.\\
This use seems to be quite elegant and synthetic, but it can create some confusions when the predicate uses the equality relation as the following
$$M_f=\underset {x\in X}E\enspace f(x)=0\enspace\footnote{page 1 of \cite{BBZ}}$$

On the definitory side, we want to notice also that, throughout the paper, it is never used the notion of ultrafilter - which is in a way more directly connected to the concept of measure - preferring the use of the dual notion of maximal ideal (which is in this context in fact nonprincipal, complete and \textit{saturated}\footnote{see \cite{J}}).

\end{document}